\documentclass[12pt]{article}
\usepackage{graphicx,latexsym, amsmath, amsfonts}
\usepackage[dvips]{color}

\def \one{\ \hbox{I\hskip-.60em 1}}

\makeatletter \@addtoreset{equation}{section}

\newtheorem{theorem}{Theorem}[section]
\newtheorem{lemma}[theorem]{Lemma}
\newtheorem{corollary}[theorem]{Corollary}
\newtheorem{prop}[theorem]{Proposition}
\newtheorem{defn}[theorem]{Definition}

\newenvironment{proof}[1][Proof]{\begin{trivlist}
\item[\hskip \labelsep {\bfseries #1}]}{\end{trivlist}}

\newtheorem{preremark}[theorem]{Remark}
  \newenvironment{remark}%
    {\begin{preremark}\upshape}{\end{preremark}}

\def\sgn{\mathop{\mathrm sgn}}

\newcommand{\eproof}{\indent\vrule height6pt width4pt depth1pt\hfil\par\medbreak}

\begin{document}

\title{\sc Necessary and sufficient conditions for the existence of the $q$-optimal measure}

\author
{ {\sc Sotirios Sabanis} \footnote{The
author is grateful to Alexander Davie and Istvan Gyongy for valuable comments and suggestions.} \\
School of Mathematics \\
University of Edinburgh, Edinburgh EH9 3JZ, UK }

\date{}

\maketitle

\begin{abstract}
This paper presents the general form and essential properties of the
$q$-optimal measure following the approach of Delbaen \&
Schachermayer (1996) and proves its existence under mild conditions.
Most importantly, it states a necessary and sufficient condition for
a candidate measure to be the $q$-optimal measure in the case even
of signed measures. Finally, an updated characterization of the
$q$-optimal measure for continuous asset price processes is
presented in the light of the counterexample appearing in Cerny \&
Kallsen (2006) concerning Hobson's (2004) approach.

\medskip \noindent

\footnotesize{ Keywords:  $q$-optimal martingale measure, uniformly
integrable martingale, signed local martingale measures, incomplete
markets}.

\end{abstract}

\section{Introduction}

In an incomplete market, the choice of the equivalent martingale
measure (EMM) for the underlying price process is not unique. Over
the last twenty years, many authors have proposed different
preference based criteria in order to choose a `suitable' pricing
measure from the class of EMMs. Two of the most popular choices are
the minimal entropy EMM, see for example Frittelli (2000), and the
variance optimal EMM, see Delbaen \& Schachermayer (1996) and
Schweizer (1996).

Recently, Hobson (2004) proposed a characterisation of the
$q$-optimal measure, for a wide range of choices of EMMs, which
includes the two aforementioned measures. The notion of
$q$-optimality is linked to the unique EMM with minimal $q$-moment
(if $q > 1$) or minimal relative entropy (if $q=1$). Hobson's (2004)
approach to identifying the $q$-optimal measure (through a so-called
fundamental equation) suggests a relaxation of an essential
condition appearing in Delbaen \& Schachermayer (1996). This
condition states that for the case $q=2$, the Radon-Nikodym process,
whose last element is the density of the candidate measure, is a
uniformly integrable martingale with respect to any EMM with a
bounded second moment. Hobson (2004) alleges that it suffices to
show that the above is true only with respect to the candidate
measure itself and extrapolates for the case $q>1$. Cerny \& Kallsen
(2006) however presented a counterexample (for $q=2$) which
demonstrates that the above relaxation does not hold in general. The
case $q=1$ is covered by Grandits \& Rheinl$\ddot{\mbox{a}}$nder
(2002).

This paper follows the approach of Delbaen \& Schachermayer (1996)
to describe and present the essential properties of the $q$-optimal
measure (with $q>1$) by extending the definition to include also
signed local martingale measures, see for example Grandits \&
Rheinl$\ddot{\mbox{a}}$nder (2002). In the light of the
counterexample appearing in Cerny \& Kallsen (2006), the analogous
sufficient condition for $q>1$ is presented to guarantee that a
candidate measure is indeed the $q$-optimal measure. Most
importantly, it is proven here that the condition under
consideration is also necessary for the identification of the
$q$-optimal measure. Furthermore, the information concerning the
form of the $q$-optimal measure helps us identify the constant
appearing in the so-called fundamental representation equation, see
Hobson (2004), which determines when a candidate measure has the
$q$-optimality property and an updated characterization of the
$q$-optimal measure is given.

\section{Main Result}

Let us consider an $\mathbb{R}^d$-valued, locally bounded, cadlag
semimartingale $S:=\{S_t\}_{t\ge 0}$ defined on a filtered
probability space $(\Omega$, $\mathcal{F}$, $\{\mathcal{F}_t\}_{t\ge
0}$, $\mathbb{P})$. It is assumed that $S$ models the evolution of
$d$ discounted stock price processes. Furthermore, let us consider
$K_0$, a linear subspace of $L^{\infty}(\mathbb{P})$, which is
spanned by simple stochastic integrals of the form (dot product)
$$
h = \phi (S_{\tau_2}-S_{\tau_1})
$$
where $\tau_1$ and $\tau_2$ are stopping times such that: (i)
$\tau_1 \le \tau_2$ a.s., (ii) the stopped process
$S^{\tau_2}:=\{S_{\tau_2\wedge t}\}_{t\ge0}$ is bounded. Moreover,
$\phi$ is assumed to be a bounded $\mathbb{R}^d$-valued
$\mathcal{F}_{\tau_1}$-measurable function. Then, we remind
ourselves of the following well-known definitions:

{\defn A probability measure $\mathbb{Q}$ on $\mathcal{F}$ with
density $u:=\frac{d\mathbb{Q}}{d\mathbb{P}}\in L^1(\mathbb{P})$ is a
local martingale measure for $S$ iff $\mathbb{Q}$ vanishes on $K_0$
i.e., $\mathbb{E}[uh] =0$} for all $h \in K_0$.

{\defn The following collection of random variables
$$
\mathcal{M}^s(\mathbb{P}) =\{u \in L^1(\mathbb{P}): \mathbb{E}[uh]=0
\mbox{ for any } h\in K_0, \mbox{ and }\mathbb{E}[u]=1\}
$$
is called the set of signed local martingale measures for the
process $S$. }

Moreover, the set of absolutely continuous (resp. equivalent) local
martingale measures $\mathcal{M}(\mathbb{P})$ (resp.
$\mathcal{M}^e(\mathbb{P})$) for the process $S$ is defined as the
intersection of $\mathcal{M}^s(\mathbb{P})$ with the positive (resp.
strictly positive) orthant of $L^1(\mathbb{P})$. Recall also here
that $\mathcal{M}^s(\mathbb{P})\cap L^q(\mathbb{P})$ is closed in
$L^q(\mathbb{P})$ and that it has a unique element of minimal
$L^q(\mathbb{P})$-norm (provided that $\mathcal{M}^s(\mathbb{P})\cap
L^q(\mathbb{P})\neq \emptyset$) due to the strict convexity of the
norm.

{\defn  Suppose that $\mathcal{M}^s(\mathbb{P})\cap
L^q(\mathbb{P})\neq \emptyset$ and $q>1$. Then, the unique element
of $\mathcal{M}^s(\mathbb{P})$ with minimal $L^q(\mathbb{P})$-norm
is called the $q$-optimal signed local martingale measure for the
process $S$}.

One can then identify the general form of the $q$-optimal measure
following the approach of Delbaen \& Schachermayer (1996). Although
this result is known in the literature, see for example Grandits
(1999), it is important in the author's view to present a relevant
proof here so as to be able to proceed with the construction of the
necessary and sufficient condition for the existence of the
$q$-optimal measure in the general framework of signed measures.

It is noted though that for $q\neq 2$, one operates in Banach spaces
instead of Hilbert spaces since the dual of $L^q(\mathbb{P})$ is
$L^p(\mathbb{P})$, where $p = \frac{q}{q-1}$. Nevertheless, it is
possible to extend Delbaen \& Schachermayer (1996) results with a
careful approach. Let $\bar{K}_0$ denote the closure of $K_0$ in
$L^p(\mathbb{P})$ and $\bar{K}$ denote the closure of the span of
$K_0$ and the constants also in $L^p(\mathbb{P})$. Then, the
annihilator of $\bar{K}_0$, which is denoted by
$\bar{K}_0^{\alpha}$, is in $L^q(\mathbb{P})$. Let also
$||\cdot||_p$ and $||\cdot||_q$ denote the $L^p(\mathbb{P})$-norm
and $L^q(\mathbb{P})$-norm respectively.

{\theorem \label{optimal} Fix $q>1$. The following statements
hold:

(a) $\mathcal{M}^s(\mathbb{P})\cap L^q(\mathbb{P})\neq \emptyset$
iff $\bar{K}_0$ does not contain the constant function 1.

(b) If $\mathcal{M}^s(\mathbb{P})\cap L^q(\mathbb{P})\neq
\emptyset$, then the probability measure $\mathbb{Q}^*$ defined by
$$\frac{d\mathbb{Q}^*}{d\mathbb{P}}:=\frac{g^*}{\mathbb{E}[g^*]},$$
where $g^*: = \sgn(1-f)|1-f|^{p-1}$ and $f$ is the unique element
of $\bar{K}_0$ with the property $$||1-f||_p
=\inf_{h\in\bar{K}_0}||1-h||_p,$$ is the unique element of
$\bar{K}_0^{\alpha}$ with minimal $L^q(\mathbb{P})$-norm.}
\begin{proof}
(a) The linear functional $\varphi \in \bar{K}_0^{\alpha}$ with
$\varphi(1) =1$ is well defined and continuous on $\bar{K}$ iff
$1\notin \bar{K}_0$.

(b) Let $f$ be the unique element of $\bar{K}_0$ such that
$||1-f||_p =\inf_{h\in\bar{K}_0}||1-h||_p$ (uniqueness is due to the
strict convexity of the $L^p(\mathbb{P})$-norm). Let $g:=1-f$, and
observe that for any other $h\in \bar{K}_0$ and $t\in \mathbb{R}$
$$
||g+th||_p^p \ge ||g||_p^p
$$
holds. As a result, we obtain
$$
\frac{d}{dt} ||g+th||_p^p |_{t=0} =0 \qquad \Rightarrow \qquad
p\mathbb{E}[\sgn(g)|g|^{p-1}h]=0.
$$
Set $g^*=\sgn(g)|g|^{p-1}$ and observe that
$\mathbb{E}[g^*]=\mathbb{E}[g^*(1-f)]=||g||_p^p>0$. Thus,
$\frac{g^*}{\mathbb{E}[g^*]} \in \bar{K}_0^{\alpha}$ and
$\mathbb{E}[\frac{g^*}{\mathbb{E}[g^*]}]=1$. Furthermore, we
calculate
$$
||\frac{g^*}{\mathbb{E}[g^*]}||_q^q = \frac{1}{||g||_p^{pq}}
\mathbb{E}[|g|^{q(p-1)}]= \frac{1}{||g||_p^q} <\infty.
$$
which implies that $\mathbb{Q}^* \in \mathcal{M}^s(\mathbb{P})\cap
L^q(\mathbb{P})$. Finally, for any element $u\in \bar{K}_0^{\alpha}$
with $\mathbb{E}[u]=1$ (i.e., any signed local martingale measure
with density in $L^q(\mathbb{P})$) we obtain
$$
\mathbb{E}[ug]= \mathbb{E}[u(1-f)]=1
$$
and thus H$\ddot{\mbox{o}}$lder inequality yields
$$
1 \le ||u||_q||g||_p\quad \Rightarrow \quad ||u||_q \ge
\frac{1}{||g||_p}=||\frac{g^*}{\mathbb{E}[g^*]}||_q.
$$
and that concludes the proof. \hfill \eproof
\end{proof}

It is the general form of the $q$-optimal measure presented in
Theorem \ref{optimal} that holds the key to obtaining the necessary
and sufficient condition for proving the $q$-optimality property of
a candidate measure. It is therefore important to recall here the
counterexample from Cerny \& Kallsen (2006). The counterexample
shows that (for $q=2$) a candidate measure may not be the
$q$-optimal measure if we only prove that the Radon-Nikodym process,
whose last element is the density of the candidate measure, is a
uniformly integrable martingale with respect to the candidate
measure itself. Therefore, we still require the condition set by
Delbaen \& Schachermayer (1996), i.e. the corresponding
Radon-Nikodym process should be a uniformly integrable martingale
with respect to any EMM with a bounded second moment. The main
Theorem of this section follows.

{\theorem \label{condition} Let $q>1$ and suppose that there exists
$\mathbb{Q}^*\in \mathcal{M}^{s}(\mathbb{P})\cap L^q(\mathbb{P})$
defined by
$$
\frac{d\mathbb{Q}^*}{d\mathbb{P}}:=\frac{g^*}{\mathbb{E}[g^*]}
$$
The following statements hold:
\begin{enumerate}
\item[(i)] if $\mathbb{Q}^*$ is the q-optimal measure, then
$\mathbb{E}_{\mathbb{Q}}[\sgn(g^*)|g^*|^{q-1}]=1$ for every
$\mathbb{Q}\in \mathcal{M}^{s}(\mathbb{P})\cap L^q(\mathbb{P})$;
\item[(ii)] conversely, if
$\mathbb{E}_{\mathbb{Q}}[\sgn(g^*)|g^*|^{q-1}]=1$ for every
$\mathbb{Q}\in \mathcal{M}^{s}(\mathbb{P})\cap L^q(\mathbb{P})$,
then $\mathbb{Q}^*$ is the $q$-optimal martingale measure.
\end{enumerate} }

\begin{proof}
(i) If $\mathbb{Q}^*$ is the q-optimal measure, then Theorem
\ref{optimal} asserts that
$$
g^* = \sgn(1-f)|1-f|^{p-1},
$$
where $||1-f||_p=\inf_{h\in \bar{K}_0}||1-h||_p$, and thus
$$
\mathbb{E}[u \sgn(g^*)|g^*|^{q-1}]= \mathbb{E}[u
\sgn(1-f)|1-f|^{(p-1)(q-1)}] = \mathbb{E}[u(1-f)]=1,
$$
for any $u \in \bar{K}_0^{\alpha}$ with $\mathbb{E}[u]=1$.

(ii) If $\mathbb{E}_{\mathbb{Q}}[\sgn(g^*)|g^*|^{q-1}]=1$ for
every $\mathbb{Q}\in \mathcal{M}^{s}(\mathbb{P})\cap
L^q(\mathbb{P})$, then
$$
\mathbb{E}[\frac{g^*}{\mathbb{E}[g^*]}\sgn(g^*)|g^*|^{q-1}]= 1
\quad \implies \mathbb{E}[g^*] = \mathbb{E}[|g^*|^q]>0.
$$
Set $\mu:= \mathbb{E}[g^*] = \mathbb{E}[|g^*|^q]$ and observe that
$$
\mathbb{E}[|\frac{d\mathbb{Q}^*}{d\mathbb{P}}|^q] =
\mathbb{E}[|\frac{g^*}{\mathbb{E}[g^*]}|^q]=\frac{\mu}{\mu^q}=\mu^{1-q}=\mu^{-q/p}.
$$
Moreover, for any $\mathbb{Q}\in \mathcal{M}^{s}(\mathbb{P})\cap
L^q(\mathbb{P})$,
$$
1 = \mathbb{E}[\frac{d\mathbb{Q}}{d\mathbb{P}}
\sgn(g^*)|g^*|^{q-1}] \le ||\frac{d\mathbb{Q}}{d\mathbb{P}}||_q
||\sgn(g^*)|g^*|^{q-1}||_p = ||\frac{d\mathbb{Q}}{d\mathbb{P}}||_q
(\mathbb{E}[|g^*|^q])^{1/p}
$$
and thus
$$
||\frac{d\mathbb{Q}}{d\mathbb{P}}||_q \ge \mu^{-1/p} \qquad
\implies \qquad \mathbb{E}[|\frac{d\mathbb{Q}}{d\mathbb{P}}|^q]
\ge \mu^{-q/p} = \mathbb{E}[|\frac{d\mathbb{Q}^*}{d\mathbb{P}}|^q]
$$
and that concludes the proof.\hfill \eproof

\end{proof}

\begin{remark} The condition $\mathbb{E}_{\mathbb{Q}}[\sgn(g^*)|g^*|^{q-1}]=1$,
which is translated as $\mathbb{E}_{\mathbb{Q}}[(g^*)^{q-1}]=1$
for the EMMs case, implies that the stochastic process
$(\hat{V}^{\mbox{\tiny opt}})^{q-1} =
\{\mathbb{E}_{\mathbb{Q}}[(V_{\infty}^{\mbox{\tiny
opt}})^{q-1}|\mathcal{F}_t]\}_{0\le t \le \infty}$ is a uniformly
integrable $\mathbb{Q}$-martingale with respect to any
$\mathbb{Q}\in \mathcal{M}^{e}(\mathbb{P})\cap L^q(\mathbb{P})$,
see Lemma \ref{uniform}. Moreover, for $q=2$, one obtains that the
corresponding Radon-Nikodym process, whose last element is the
density $\frac{d\mathbb{Q}^*}{d\mathbb{P}}$, is a uniformly
integrable $\mathbb{Q}$-martingale with respect to any
$\mathbb{Q}\in \mathcal{M}^{e}(\mathbb{P})\cap L^2(\mathbb{P})$
and this is a necessary and sufficient condition for
$\mathbb{Q}^*$ to be the $q$-optimal (local) martingale measure.
\end{remark}

Let us turn our attention now to the case where $S$ is a continuous
adapted stochastic process. Then, one can prove that $\mathbb{Q}^*$
is a probability measure equivalent to $\mathbb{P}$. This result is
also known in the literature, see for example  Grandits \&
Rheinl$\ddot{\mbox{a}}$nder (2002), but it is presented here as the
generalisation of Delbaen \& Schachermayer (1996) technique for
completeness of this section.

Note also that the notation $(\varphi\cdot S)_t\in \bar{K}_0$ is
used as a shorthand notation for the stochastic integral
$$
(\varphi\cdot S)_t = \int_0^t \varphi_u dS_u
$$
for every $0\le t\le \infty$, where the process $\varphi \in
\mathcal{H}_p$, i.e. it satisfies
$$
\mathbb{E}[(\int_0^{\infty}\varphi^2_t d[S]_t)^{p/2}]<\infty.
$$

{\theorem \label{non-neg} Fix $q>1$. Let us assume that $S$ is a
continuous process and that $\mathcal{M}^s(\mathbb{P})\cap
L^q(\mathbb{P})\neq \emptyset$. Then, the $q$-optimal signed local
martingale measure $\mathbb{Q}^*$ is a well-defined probability
measure absolutely continuous with $\mathbb{P}$}.
\begin{proof}
In order to show that $\frac{d\mathbb{Q}^*}{d\mathbb{P}}$ is
non-negative, it suffices to prove that $f\le1$ (a.s.).

Let us assume (on the contrary) that there exists $\epsilon \in (0,
1)$ such that $\mathbb{P}(f>1+\epsilon)>\epsilon$. Then, there
exists a simple integrand $\phi$ such that
\begin{enumerate}
\item[(a)] $(\phi\cdot S)_{\infty} \in \bar{K}_0$,
\item[(b)] $||(\phi\cdot S)_{\infty}-f||_p \le c$, where $c<
(\frac{\epsilon}{2})^{p+1}(\sum_{i=1}^{\lceil p/2 \rceil}
\binom{p}{2i-1})^{-1}$ and,
\item[(c)] $||1-(\phi\cdot S)_{\infty}||_p\le 1$ (since $||1-f||_p\le 1<1+c$).
\end{enumerate}
Then, we observe that
$$
\mathbb{P}((\phi\cdot S)_{\infty}>1+ \frac{\epsilon}{2}) \ge
\mathbb{P}(f>1+\epsilon) - \mathbb{P}(|f-(\phi\cdot
S)_{\infty}|>\frac{\epsilon}{2})
$$
and since
$$
\mathbb{P}(|f-(\phi\cdot
S)_{\infty}|>\frac{\epsilon}{2}) \le
(\frac{2}{\epsilon})^p\mathbb{E}[|f-(\phi\cdot S)_{\infty}|^p] \le
(\frac{2}{\epsilon})^p c^p
$$
we conclude that
$$
\mathbb{P}((\phi\cdot S)_{\infty}>1+ \frac{\epsilon}{2}) \ge
\epsilon -
(\frac{2}{\epsilon})^p(\frac{\epsilon}{2})^{(p+1)p}(\sum_{i=1}^{\lceil
p/2 \rceil}\binom{p}{2i-1})^{-p} \ge \frac{\epsilon}{2}.
$$
Moreover, we define the stopping time $\tau=\inf\{t\ge0: (\phi\cdot
S)_t>1\}$. Then,
$$
|1-(\phi\cdot S)_{\infty}|^p = |1-(\phi\cdot
S)_{\tau}|^p\one_{\{\tau=\infty\}} + |1-(\phi\cdot
S)_{\infty}|^p\one_{\{\tau<\infty\}} = |1-(\phi\cdot S)_{\tau}|^p +
|1-(\phi\cdot S)_{\infty}|^p\one_{\{\tau<\infty\}}
$$
since for $\tau <\infty$ we have $1-(\phi\cdot S)_{\tau}=0$ due to
the continuity of $S$. Hence,
\begin{align}
||1-(\phi\cdot S)_{\infty}||_p^p & = ||1-(\phi\cdot S)_{\tau}||_p^p
+ \mathbb{E}[|1-(\phi\cdot S)_{\infty}|^p\one_{\{\tau<\infty\}}]
\nonumber \\ & \ge  ||1-(\phi\cdot S)_{\tau}||_p^p +
\mathbb{E}[|1-(\phi\cdot S)_{\infty}|^p\one_{\{(\phi\cdot
S)_{\infty}>1+ \frac{\epsilon}{2}\}}]\nonumber \\ & \ge
||1-(\phi\cdot S)_{\tau}||_p^p +
(\frac{\epsilon}{2})^p\mathbb{P}((\phi\cdot S)_{\infty}>1+
\frac{\epsilon}{2}) \nonumber \\ &  \ge ||1-(\phi\cdot
S)_{\tau}||_p^p + (\frac{\epsilon}{2})^{p+1}. \nonumber
\end{align}
Note also that due to Minkowski inequality
$$
||1-(\phi\cdot S)_{\infty}||_p \le ||1-f||_p + ||(\phi\cdot
S)_{\infty}-f||_p \le ||1-f||_p + c
$$
which implies
\begin{align}
||1-f||_p^p &\ge ||1-(\phi\cdot S)_{\infty}||_p^p +
\sum_{i=1}^p\binom{p}{i}(-1)^i||1-(\phi\cdot
S)_{\infty}||_p^{p-i}c^i \nonumber \\ & \ge ||1-(\phi\cdot
S)_{\infty}||_p^p - \sum_{i=1}^{\lceil p/2
\rceil}\binom{p}{2i-1}c^{2i-1} \nonumber \\ & \ge ||1-(\phi\cdot
S)_{\tau}||_p^p + (\frac{\epsilon}{2})^{p+1} - \sum_{i=1}^{\lceil
p/2 \rceil}\binom{p}{2i-1}c^{2i-1} \nonumber \\ & \ge
||1-(\phi\cdot S)_{\tau}||_p^p + (\frac{\epsilon}{2})^{p+1} -
c\sum_{i=1}^{\lceil p/2 \rceil}\binom{p}{2i-1} \nonumber \\ & >
||1-(\phi\cdot S)_{\tau}||_p^p \nonumber
\end{align}
which is a contradiction since $f$ is the unique element of
$\bar{K}_0$ with the property $||1-f||_p
=\inf_{h\in\bar{K}_0}||1-h||_p$. \hfill \eproof
\end{proof}

Theorem \ref{optimal} states that $f\in \bar{K}_0$, therefore under
the assumption that $S$ is a semi-martingale, we can represent
$$
g = 1-f =1-(\psi\cdot S)_{\infty}.
$$
Moreover, we fix $\mathbb{Q}\in \mathcal{M}^e(\mathbb{P})\cap
L^q(\mathbb{P})$ and for every $t\ge 0$ we define
\begin{align}
V_{\infty}^{opt}:= \frac{g^*}{\mathbb{E}[g^*]} =
\frac{g^{p-1}}{\mathbb{E}[g^*]} \quad &\& \quad V_t^{opt}=
\mathbb{E}[V_{\infty}^{opt}|\mathcal{F}_t],\nonumber
\\ X_t := \mathbb{E}[g^{p-1}|\mathcal{F}_t] =
\mathbb{E}[g^*]V_t^{opt} \quad &\& \quad Y_t = 1-(\psi\cdot S)_t =
\mathbb{E}_{\mathbb{Q}}[g|\mathcal{F}_t] \nonumber
\end{align}
and the stopping times
$$
\tau :=\inf\{t\ge0: X_t=0\}  \quad \& \quad \sigma:=\inf\{t\ge0:
Y_t=0\}.
$$
Note that the processes $X$ and $Y$ are non-negative
supermartingales with non-negative last elements $X_{\infty}$ and
$Y_{\infty}$, therefore when any of their paths hits zero, it stays
at zero. Furthermore, the continuity of $Y$ implies that the
stopping time $\sigma$ is predictable. As a result, the following
lemmas (\ref{1} and \ref{folkore}) can be proved in a similar
fashion as in Delbaen \& Schachermayer (1996).

{\lemma \label{1} Fix $q>1$. Let us assume that $S$ is a continuous
semi-martingale and that $\mathcal{M}^e(\mathbb{P})\cap
L^q(\mathbb{P})\neq \emptyset$. Then, $\tau = \sigma$.}

\begin{proof}
Our aim is to prove that $\mathbb{P}(\sigma<\tau) =
\mathbb{P}(\sigma>\tau)=0$. Consider the set $\{\sigma<\tau\}$, then
$$
0<X_{\sigma}= \mathbb{E}[X_{\infty}|\mathcal{F}_{\sigma}] =
\mathbb{E}[g^{p-1}|\mathcal{F}_{\sigma}] = \mathbb{E}[(1-(\psi\cdot
S)_{\sigma})^{p-1}|\mathcal{F}_{\sigma}]
$$
since $1-(\psi\cdot S)_{\infty} = 1-(\psi\cdot S)_{\sigma}$ on
$\{\sigma<\tau\} \subset \{\sigma<\infty\}$ and thus
$$
0<X_{\sigma}=(1-(\psi\cdot S)_{\sigma})^{p-1}=0 \qquad
(\mbox{contradiction})
$$
which implies $\mathbb{P}(\sigma<\tau) =0$. Now consider the set
$\{\sigma>\tau\}\subset \{\tau<\infty\}$ and observe that
$$
0=X_{\tau} = \mathbb{E}[X_{\infty}|\mathcal{F}_{\tau}]=
\mathbb{E}[g^{p-1}|\mathcal{F}_{\tau}]
$$
which implies $g =0$ on $\{\sigma>\tau\}$. Thus, since $Y_{\tau} =
\mathbb{E}_{\mathbb{Q}}[g|\mathcal{F}_{\tau}]$, we obtain $Y_{\tau}
=0$ on $\{\sigma>\tau\}$ (contradiction) which implies
$\mathbb{P}(\sigma>\tau) =0$. \hfill \eproof
\end{proof}

{\corollary The martingale $V^{opt}$ is continuous at time $t=\tau$
and the stopping time $\tau$ is predictable and thus is announced by
the sequence $\tau_n=\inf\{t\ge 0: V^{opt}_t \le \frac{1}{n}\}\wedge
n$. }

{\lemma \label{folkore} Let $M:=\{M_t\}_{0\le t\le \infty}$ be a
$q$th integrable martingale such that $M_0>0$ . Let also
$\tau=\inf\{t\ge 0: M_t =0\}$ be a predictable stopping time
announced by a sequence of stopping times $\{\tau_n\}_{n\ge1}$.
Then,
$$
\mathbb{E}[\frac{M^q_{\infty}}{M^q_{\tau_n}}|\mathcal{F}_{\tau_n}]
\to \infty
$$
on the set $\{M_{\tau}=0\}$.}
\begin{proof}
First observe that
$$
\one=\mathbb{E}[\frac{M_{\infty}}{M_{\tau_n}}|\mathcal{F}_{\tau_n}]
= \mathbb{E}[\frac{M_{\infty}}{M_{\tau_n}}\one_{\{M_{\tau}\neq
0\}}|\mathcal{F}_{\tau_n}] \le
\mathbb{E}[(\frac{M_{\infty}}{M_{\tau_n}})^q|\mathcal{F}_{\tau_n}]^{1/q}
\mathbb{E}[\one_{\{M_{\tau}\neq 0\}}|\mathcal{F}_{\tau_n}]^{1/p}
$$
and then recall that $\mathbb{E}[\one_{\{M_{\tau}\neq
0\}}|\mathcal{F}_{\tau_n}]$ tends to zero on $\{M_{\tau}=0\}$.
\hfill \eproof
\end{proof}

{\theorem  Fix $q>1$. Let us assume that $S$ is a continuous
semi-martingale and that $\mathcal{M}^e(\mathbb{P})\cap
L^q(\mathbb{P})\neq \emptyset$. Then, the $q$-optimal local
martingale measure $\mathbb{Q}^*$ is in fact equivalent to
$\mathbb{P}$.}
\begin{proof}
Let us assume on the contrary that $\mathbb{P}[X_{\tau}=0]>0$ and
observe that for the uniformly integrable martingale $V$, where
$V_t:=\mathbb{E}[\frac{d\mathbb{Q}}{d\mathbb{P}}|\mathcal{F}_t]$
for all $0\le t\le \infty$ and $\mathbb{Q} \in
\mathcal{M}^e(\mathbb{P})\cap L^q(\mathbb{P})$, we have
$\inf_{t\ge0}V_t>0$ and $\sup_{t\ge
0}\mathbb{E}[(V_{\infty})^q|\mathcal{F}_t]<\infty$ (both
inequalities hold a.s.). In view of Lemma \ref{folkore}, one
expects that for a large enough $n$ the set
$$
A=\Big\{\sup_{t\ge0}\frac{\mathbb{E}[(V_{\infty})^q|\mathcal{F}_t]}{(V_t)^q}
< \frac{\mathbb{E}[(V_{\infty}^{\mbox{\tiny
opt}})^q|\mathcal{F}_{\tau_n}]}{(V_{\tau_n}^{\mbox{\tiny
opt}})^q}\Big\}
$$
is non empty, thus
$$
A_n=\Big\{\frac{\mathbb{E}[(V_{\infty})^q|\mathcal{F}_{\tau_n}]}{(V_{\tau_n})^q}
< \frac{\mathbb{E}[(V_{\infty}^{\mbox{\tiny
opt}})^q|\mathcal{F}_{\tau_n}]}{(V_{\tau_n}^{\mbox{\tiny
opt}})^q}\Big\}
$$
is non empty in $\mathcal{F}_{\tau_n}$. Then, the martingale
\begin{align}
\bar{V}_t=
\begin{cases}
V_t^{\mbox{\tiny opt}}, & t<\tau_n, \\
\frac{V_t V_{\tau_n}^{\mbox{\tiny opt}}}{V_{\tau_n}}, &
\mbox{for } t\ge\tau_n \mbox{ on the set }A_n, \\
V_t^{\mbox{\tiny opt}}, & \mbox{for } t\ge\tau_n \mbox{ on the
complement of the set }A_n,
\end{cases}
\nonumber
\end{align}
defines an equivalent martingale measure $\bar{\mathbb{Q}}$ to
$\mathbb{P}$ such that $||\bar{V}_{\infty}||_q <
||V_{\infty}^{\mbox{\tiny opt}}||_q$ which is clearly a
contradiction. \hfill \eproof
\end{proof}

The last Lemma of this section provides the connection between the
condition appearing in Theorem \ref{condition} and the behaviour of
$\{(\hat{V}_t^{\mbox{\tiny opt}})^{q-1}\}_{0\le t\le\infty}$ as
defined below.

{\lemma \label {uniform} Fix $q>1$. Let
$\mathcal{M}^{e}(\mathbb{P})\cap L^q(\mathbb{P})\neq \emptyset$
and fix $\hat{\mathbb{Q}}\in \mathcal{M}^{e}(\mathbb{P})\cap
L^q(\mathbb{P})$. Let us define the process $\hat{V}$ by
$\hat{V}_t^{\mbox{\tiny
opt}}:=(\mathbb{E}_{\hat{\mathbb{Q}}}[(V_{\infty}^{\mbox{\tiny
opt}})^{q-1}|\mathcal{F}_t])^{1/(q-1)}$ for every $t\ge 0$. Then,
\begin{equation}
(\hat{V}_t^{\mbox{\tiny opt}})^{q-1} = ||V_{\infty}^{\mbox{\tiny
opt}} ||_q^q+ (\varphi\cdot S)_t
\end{equation}
where the stochastic integral $(\varphi\cdot S)$ is well defined,
i.e. $\varphi \in \mathcal{H}_p$, and is a uniformly integrable
$\mathbb{Q}$-martingale  for every $\mathbb{Q}\in
\mathcal{M}^{e}(\mathbb{P})\cap L^q(\mathbb{P})$. Furthermore, the
choice of $\varphi$ is independent of the choice of
$\hat{\mathbb{Q}}\in \mathcal{M}^{e}(\mathbb{P})\cap
L^q(\mathbb{P})$. }

\begin{proof}
Recall that $g\in \bar{K}$ and $(g^*)^{q-1}=g$ which imply that
there exists a sequence $\{g_i\}_{i\ge1} \in K$ that converges to
$(V_t^{\mbox{\tiny opt}})^{q-1}$ in $L^p(\mathbb{P})$. Moreover, we
observe that
$$
\mathbb{E}_{\hat{\mathbb{Q}}}[g_i - (V_{\infty}^{\mbox{\tiny
opt}})^{q-1}] = \mathbb{E}[(g_i - (V_{\infty}^{\mbox{\tiny
opt}})^{q-1})\frac{d\hat{\mathbb{Q}}}{d\mathbb{P}}] \le ||g_i -
(V_{\infty}^{\mbox{\tiny opt}})^{q-1}
||_p||\frac{d\hat{\mathbb{Q}}}{d\mathbb{P}} ||_q
$$
which implies convergence in $L^1(\hat{\mathbb{Q}})$. Note that if
we choose to represent each $g_i \in K$ as follows
$$
g_i = \delta_i + (\phi_i\cdot S)
$$
where $\delta_i$ denotes the real number in the representation, we
obtain as a result that
\begin{align}
\lim_{i\to\infty} \delta_i & = \lim_{i\to\infty}
\mathbb{E}_{\hat{\mathbb{Q}}}[g_i] = \mathbb{E}_{\hat{\mathbb{Q}}}[
(V_{\infty}^{\mbox{\tiny opt}})^{q-1}] = \mathbb{E}[
(V_{\infty}^{\mbox{\tiny
opt}})^{q-1}\frac{d\hat{\mathbb{Q}}}{d\mathbb{P}}] = \frac{1}{(\mathbb{E}[g^*])^{q-1}} \nonumber \\
& = \mathbb{E}[ (V_{\infty}^{\mbox{\tiny
opt}})^{q-1}\frac{d\mathbb{Q}^*}{d\mathbb{P}}] =
\frac{1}{||g||_p^q} = ||V_{\infty}^{\mbox{\tiny opt}} ||_q^q,
\nonumber
\end{align}
so the process $\{g_i-\delta_i\}_{1\le i \le \infty}$ converges in
$L^1(\hat{\mathbb{Q}})$ to $(V_{\infty}^{\mbox{\tiny opt}})^{q-1} -
||V_{\infty}^{\mbox{\tiny opt}} ||_q^q$. Thus, following once more
the approach of Delbaen \& Schachermayer (1996), one obtains that
the choice of $\varphi$ is independent of the choice of
$\hat{\mathbb{Q}}$ since the process $(\varphi\cdot S)$ is a
uniformly integrable $\mathbb{Q}$-martingale  for every
$\mathbb{Q}\in \mathcal{M}^{e}(\mathbb{P})\cap L^q(\mathbb{P})$
converging to $(V_{\infty}^{\mbox{\tiny opt}})^{q-1} -
||V_{\infty}^{\mbox{\tiny opt}} ||_q^q$ in $L^1(\mathbb{Q})$. \hfill
\eproof
\end{proof}

\section{Continuous Univariate Case}

Let $T\in (0,\infty]$ denote the termination date of the economy,
i.e. we can work under either a finite ($T<\infty$) or an infinite
($T=\infty$) time horizon. Let $(\Omega$, $\mathcal{F}$,
$\{\mathcal{F}_t\}_{0 \le t\le T}$, $\mathbb{P})$ be a filtered
probability space that satisfies the usual conditions of
right-continuity and completeness, where
$\mathcal{F}=\mathcal{F}_T$ and $\mathcal{F}_0$ is trivial.
Moreover, let $Y:=\{Y_t\}_{0 \le t\le T}$ denote the volatility of
the traded asset $S$. Suppose that $S$ is a continuous
semimartingale governed by the following stochastic differential
equation
\begin{equation}
\label{main} \qquad \qquad dS_t  = \mu(S_t,Y_t,t)  dt +
\sigma(S_t, Y_t,t)dB_t, \qquad \forall \mbox{ } t\in [0,T],
\end{equation}
where $B:=\{B_t\}_{0 \le t\le T}$ is a $\mathbb{P}$-Brownian
motion. The semimartingale $S$ admits a Doob-Meyer decomposition
given by
\begin{equation} \label{semimartingale}
S = S_0 +  A^S + M^S
\end{equation}
where $A^S$ denotes an increasing process and $M^S$ denotes a local
martingale. Furthermore, consider the processes
$$
\lambda:= \frac{\mu}{\sigma}, \qquad
\bar{\lambda}:=\frac{\lambda}{\sigma} \qquad \& \qquad
\bar{\eta}:=\frac{\eta}{\sigma}
$$
and observe that in the context of equation (\ref{main})
$$
A^S_t:= \int_0^t\mu_tdt, \qquad M^S_t = \int_0^t\sigma_tdB_t \qquad
\& \qquad A^S = \bar{\lambda} \cdot [M^S]
$$
for all $t\in [0,T]$. Then, the following proposition sets out
sufficient criteria so that a candidate measure should satisfy in
order to be the $q$-optimal measure.

{\prop \label{criteria} Let $T\in (0,\infty]$ and $q>1$ be fixed.
Suppose that there exists a $B$-integrable, predictable process
$\eta$ such that
\begin{enumerate}
\item[(i)] $\mathbb{E}_{\mathbb{Q}}[\mathcal{E}((q-1)(\bar{\eta}-\bar{\lambda})\cdot
S)_T]=1$ for every $\mathbb{Q}\in \mathcal{M}^{e}(\mathbb{P})\cap
L^q(\mathbb{P})$,
\item[(ii)] $\mathbb{E}[(\mathcal{E}((q-1)(\bar{\eta}-\bar{\lambda})\cdot
S)_T)^{p-1}]$ is a non-zero finite constant and,
\item[(ii)] it satisfies
\begin{equation} \label{fundamental}
\exp(\frac{q}{2}\bar{\lambda}\cdot A^S_T)\mathcal{E}(M^Y)_T = c
\mathcal{E}(\bar{\eta}\cdot(M^S+qA^S))_T \exp(-\frac{q-2}{2}
\bar{\eta}^2\cdot[M^S]_T),
\end{equation}
where $M^Y$ is a a local martingale with $<M^S, M^Y>=0$ and $c$ is
given by
$$c=
1/\mathbb{E}[(\mathcal{E}((q-1)(\bar{\eta}-\bar{\lambda})\cdot
S)_T)^{p-1}].$$
\end{enumerate}
Then, $V^{\mbox{\tiny opt}}:= \mathcal{E}(-\lambda\cdot B - M^Y)$
is a uniformly integrable $\mathbb{P}$-martingale, and
$\mathbb{Q}^*$ with density $V_T^{\mbox{\tiny opt}}$ is the
$q$-optimal measure.}

\begin{proof}
The integrability condition imposed on $\eta$ guarantees the
existence of the stochastic integrals appearing in equation
(\ref{fundamental}). Then, we calculate
\begin{align}
V_T^{\mbox{\tiny opt}} & = \mathcal{E}(-\lambda\cdot B - M^Y)_T =
\mathcal{E}(-\bar{\lambda}\cdot M^S)_T
\exp(-\frac{q}{2}\bar{\lambda}\cdot A^S_T)
\exp(\frac{q}{2}\bar{\lambda}\cdot A^S_T)\mathcal{E}(M^Y)_T
\nonumber \\
& = \exp(-\bar{\lambda}\cdot S_T)
\exp(-\frac{q-1}{2}\bar{\lambda}\cdot A^S_T)c
\mathcal{E}(\bar{\eta}\cdot(M^S+qA^S))_T
\exp(-\frac{q-2}{2}\bar{\eta}^2\cdot[M^S]_T) \nonumber \\
& = \exp((\bar{\eta}-\bar{\lambda})\cdot S_T)
\exp(-\frac{q-1}{2}\bar{\lambda}\cdot A^S_T)c
\exp((q-1)\bar{\eta}\bar{\lambda}\cdot[M^S]_T)
\exp(-\frac{q-1}{2}\bar{\eta}^2\cdot[M^S]_T) \nonumber \\
& = c \exp((\bar{\eta}-\bar{\lambda})\cdot S_T -
\frac{q-1}{2}(\bar{\eta}-\bar{\lambda})^2\cdot [S]_T) = c
(\mathcal{E}((q-1)(\bar{\eta}-\bar{\lambda})\cdot S)_T)^{p-1}
\nonumber
\end{align}
and consequently $\mathbb{Q}^*\in \mathcal{M}^{e}(\mathbb{P})$
since $\mathbb{E}[V_T^{\mbox{\tiny opt}}] = \mathbb{E}[c
(\mathcal{E}((q-1)(\bar{\eta}-\bar{\lambda})\cdot S)_T)^{p-1}]=1$
due to condition (ii). Moreover, $\mathbb{Q}^* \in
L^q(\mathbb{P})$ since
$$
(\frac{d\mathbb{Q}^*}{d\mathbb{P}})^{q-1} = (V_T^{\mbox{\tiny
opt}})^{q-1} =
 c^{q-1}\mathcal{E}((q-1)(\bar{\eta}-\bar{\lambda})\cdot S)_T
$$
which yields
$$
\mathbb{E}[(\frac{d\mathbb{Q}^*}{d\mathbb{P}})^q]=
\mathbb{E}_{\mathbb{Q}^*}[(\frac{d\mathbb{Q}^*}{d\mathbb{P}})^{q-1}]=
\mathbb{E}_{\mathbb{Q}^*}[c^{q-1}\mathcal{E}((q-1)(\bar{\eta}-\bar{\lambda})\cdot
S)_T] = c^{q-1} <\infty.
$$
Condition (i) and Theorem \ref{condition} assert that $\mathbb{Q}^*$
is the $q$-optimal martingale measure. Furthermore, Theorem
\ref{optimal} identifies $g$ as the last element
$\mathcal{E}((q-1)(\bar{\eta}-\bar{\lambda})\cdot S)_T$ of the
uniformly integrable $\mathbb{Q}$-martingale
$\mathcal{E}((q-1)(\bar{\eta}-\bar{\lambda})\cdot S)$. \hfill
\eproof
\end{proof}
\begin{remark}
Another byproduct of the $q$-optimal measure comes from
$$
1 =
\mathbb{E}_{\mathbb{Q}}[\mathcal{E}((q-1)(\bar{\eta}-\bar{\lambda})\cdot
S)_T]= \mathbb{E}[c
(\mathcal{E}((q-1)(\bar{\eta}-\bar{\lambda})\cdot S)_T)^{p-1}
\mathcal{E}((q-1)(\bar{\eta}-\bar{\lambda})\cdot S)_T],
$$
which yields
$$
\mathbb{E}[(\mathcal{E}((q-1)(\bar{\eta}-\bar{\lambda})\cdot
S)_T)^p]
=\mathbb{E}[(\mathcal{E}((q-1)(\bar{\eta}-\bar{\lambda})\cdot
S)_T)^{p-1}].
$$
A property that holds also due to $\mathbb{E}[g^p]
=\mathbb{E}[g^{p-1}(1-f)]=\mathbb{E}[g^{p-1}]$ according to
Theorem \ref{optimal}.
\end{remark}
\begin{remark}
Equation (\ref{fundamental}) in Proposition \ref{criteria} is a
generalisation of the Fundamental Equation (1.2) in Hobson (2004)
and Equation (3.2) in Cerny \& Kallsen. Moreover, condition {\it
(i)} in Proposition \ref{criteria} is the essential difference with
Theorem 3.1, page 543, in Hobson (2004) and addresses the issue
related to the counterexample presented by Cerny \& Kallsen (2006).
Condition {\it (i)} is replaced by the weaker condition
$$
\mathbb{E}_{\mathbb{Q}^{(q)}}[\mathcal{E}((q-1)(\bar{\eta}-\bar{\lambda})\cdot
S)_T]=1,
$$
where $\mathbb{Q}^{(q)}$ is a candidate measure, in  Hobson (2004).
\end{remark}
\begin{remark}
In Hobson (2004), $Y$ is assumed to be driven by
\begin{equation} \label{vol}
dY_t  = \alpha(Y_t, t)dt + \beta(Y_t, t)dW_t, \qquad \forall
\mbox{ } t\in[0,T],
\end{equation}
which implies that $M^Y=\xi\cdot W $, where $W:=\{W_t\}_{0 \le t\le
T}$ is a $\mathbb{P}$-Brownian motions such that $ dW_t = \rho_t
dB_t + \sqrt{1-\rho^2_t}dZ_t, $ $B$ and $Z:=\{Z _t\}_{0 \le t\le T}$
are independent $\mathbb{P}$-Brownian motions and $\rho_t$ is the
instantaneous correlation. It is possible then to identify the
constant $c_H$ appearing in Hobson's so-called fundamental
representation equation, i.e. equation (1.2), page 538,
\begin{equation} \label{constant}
c_H = \ln c =
-\ln(\mathbb{E}[(\mathcal{E}((q-1)(\bar{\eta}-\bar{\lambda})\cdot
S)_T)^{p-1}]) =
-\ln(\mathbb{E}[(\mathcal{E}((q-1)(\bar{\eta}-\bar{\lambda})\cdot
S)_T)^p])
\end{equation}
and observe that indeed
$$
\mathbb{E}[(\frac{d\mathbb{Q}^*}{d\mathbb{P}})^q] =c^{q-1} =
e^{c_H(q-1)}.
$$
Moreover, for the case where $\lambda_t=\lambda(t)$, i.e $\lambda$
is only a deterministic function of time, $\eta\equiv\xi\equiv0$ is
the solution to equation (\ref{fundamental}), and immediately one
derives that
$$
c_H = \frac{q}{2}\int_0^T\lambda^2(t)dt
$$
which is also obtained by equation (\ref{constant}) and agrees with
the findings in Hobson (2004).
\end{remark}


\begin{remark}
Similarly, let us suppose that equation (\ref{vol}) holds and
moreover, $B$ and $W$ are independent,
$\lambda_t\equiv\lambda(Y_t,t)$, i.e. $\mu(S_t,Y_t,t) =
\hat{\mu}(Y_t,t)S_t$ and $\sigma(S_t,Y_t,t) =
\hat{\sigma}(Y_t,t)S_t$, and the ``mean-variance trade-off process"
$K_t:= \int_0^t\lambda_t^2dt$ is uniformly bounded, then one obtains
the same result as in the example appearing in pages 1032--1036 in
Grandits \& Rheinl$\ddot{\mbox{a}}$nder (2002). It is an immediate
consequence of Proposition \ref{criteria}.

In order to highlight the importance of Proposition \ref{criteria}
and prove the above claim, observe that for $\eta\equiv 0$
conditions {\it (i)} and {\it (ii)} are immediately satisfied and
equation (\ref{fundamental}) is reduced to
$$
\mathcal{E}(M^Y)_T = c\exp(-\frac{q}{2}\int_0^T \lambda^2(Y_t,t)dt).
$$
Then, the Martingale Representation Theorem guarantees that there
exists a solution. As a result, all conditions of Proposition
\ref{criteria} are satisfied and
$$
V_T^{\mbox{\tiny opt}} = c (\mathcal{E}(-(q-1)\bar{\lambda}\cdot
S)_T)^{p-1}
$$
is the $q$-optimal measure. Moreover, $V_T^{\mbox{\tiny opt}}$ can
be rewritten as
$$
V_T^{\mbox{\tiny opt}} = c
\exp\Big(-\frac{1}{2}\Big(1+\frac{1}{p-1}\Big)\int_0^T\frac{\mu^2_t}{\sigma^2_t}dt\Big)
\mathcal{E}\Big(-\frac{\mu}{\sigma}\cdot W\Big)_T
$$
which is the same as the representation given in Grandits \&
Rheinl$\ddot{\mbox{a}}$nder (2002), page 1034. Futhermore, one can
show
\begin{align}
c_H & =
-\ln(\mathbb{E}[(\mathcal{E}((q-1)(\bar{\eta}-\bar{\lambda})\cdot
S)_T)^p])  =
-\ln(\mathbb{E}[(\mathcal{E}((q-1)(\bar{\eta}-\bar{\lambda})\cdot
S)_T)^{p-1}]) \nonumber \\ & =
-\ln(\mathbb{E}[\exp(-\frac{q}{2}K_T)]) \nonumber
\end{align}
which agrees with the findings in Hobson (2004).
\end{remark}


\begin{thebibliography}{9}

\bibitem{Cerny}
{\sc Cerny, A.} and {\sc J. Kallsen} (2006), A Counterexample
Concerning The Variance-Optimal Martingale Measure, {\it
http://ssrn.com/abstract=912952}, to appear in Mathematical Finance.

\bibitem{Schacher1}
{\sc Delbaen, F.,} and {\sc W. Schachermayer} (1994): A General
Version of the Fundamental Theorem of Asset Pricing, {\it Math.
Annalen}, {\bf 300}, 463--520.

\bibitem{Schacher2}
{\sc Delbaen, F.,} and {\sc W. Schachermayer} (1996):  The
Variance-Optimal Martingale Measure for Continuous Processes, {\it
Bernoulli}, {\bf 2}, 81--106.

\bibitem{Delbaen et al}
{\sc Delbaen, F., P. Monat, W. Schachermayer, M. Schweizer}, and
{\sc C. Stricker} (1997): Weighted norm inequalities and hedging in
incomplete markets, {\it Finance and Stochastics}, {\bf 1},
181--227.

\bibitem{Frittelli}
{\sc Frittelli, M.} (2000): The Minimal Entropy Measure and the
Valuation Problem in Incomplete Markets, {\it Mathematical Finance},
{\bf 10}, 39--52.

\bibitem{Grandits}
{\sc Grandits, P.} (1999): The $p$-Optimal Martingale Measure and
its Asymptotic Relation with the Minimal Entropy Martingale Measure,
{\it Bernoulli}, {\bf 5}(2), 225--247.


\bibitem{Grandits2}
{\sc Grandits, P.,} and {\sc T. Rheinlander} (2002): On the minimal
entropy martingale measure, {\it The Annals of Probability}, {\bf
30}, 1003–-1038.

\bibitem{Hobson}
{\sc Hobson, D.} (2004): Stochastic Volatility Models, Correlation,
and the $q$-Optimal Measure, {\it Mathematical Finance}, {\bf 14},
537--556.

\bibitem{Karatzas}
{\sc Karatzas, I.,} and {\sc S. Shreve} (1988): {\it Brownian Motion
and Stochastic Calculus}, New York: Springer-Verlag.

\bibitem{Schweizer}
{\sc Schweizer, M.,} (1996): Approximation Pricing and the
Variance-Optimal Martingale Measure,  {\it The Annals of
Probability}, {\bf 24}, 206--236.




\end{thebibliography}
\end{document}